\documentclass[12pt]{amsart}

\usepackage{amsmath, amscd, graphicx, latexsym, hyperref, times, rlepsf}

\oddsidemargin .25in \theoremstyle{plain}

\newtheorem{Thm}{Theorem}

\newtheorem{Cor}[Thm]{Corollary}

\newtheorem{Prop}[Thm]{Proposition}

\newtheorem{Def}[Thm]{Definition}

\newtheorem{Rem}[Thm]{Remark}

\newcommand{\genus}{\operatorname{genus}}

\newcommand{\sg}{\operatorname{sg}}

\newcommand{\sn}{\operatorname{sn}}

\newcommand{\bn}{\operatorname{bn}}

\newcommand{\Mg}{\operatorname{Mg}}

\newcommand{\Hg}{\operatorname{Hg}}

\newcommand{\Mb}{\operatorname{Mb}}

\newcommand{\Mn}{\operatorname{Mn}}

\newcommand{\bc}{\operatorname{bc}}

\newcommand{\Div}{\operatorname{Div}}

\newcommand{\rot}{\operatorname{rot}}

\newcommand{\rk}{\operatorname{rk}}

\newcommand{\tb}{\mathop{\rm tb}\nolimits}

\newcommand{\bfz}{{\mathbb{Z}}}

\newcommand{\bfc}{{\mathbb{C}}}

\newcommand{\bfq}{{\mathbb{Q}}}

\newcommand{\bfr}{{\mathbb{R}}}

\newcommand{\OB}{\mathcal{OB}}

\def\bdy{\partial}

\def\v{\vskip.12in}

\def\a{\alpha}
\def\b{\beta}
\def\d{\delta}
\def\g{\gamma}
\def\G{\Gamma}

\begin{document}

\title[]{Milnor open books of links of some rational surface singularities}

\author{Mohan Bhupal}

\author{Burak Ozbagci}

\address{Department of Mathematics \\ METU \\ Ankara \\ Turkey}

\email{bhupal@metu.edu.tr}

\address{Mathematical Sciences Research Institute \\ 17 Gauss Way \\ Berkeley, CA, 94720-5070}

\email{bozbagci@msri.org}



\begin{abstract}

We describe Milnor open books and Legendrian
surgery diagrams for canonical contact structures of links of some rational surface singularities.
We also describe an infinite  family of  Milnor fillable contact $3$-manifolds so that
the Milnor genus (resp.\ Milnor norm) is strictly greater than the support genus (resp.\ support norm) of the canonical contact structure, for each member of this family.

\end{abstract}

\maketitle

\section{Introduction}

The link of a normal complex surface singularity  carries a canonical contact structure $\xi_{can}$ (a.k.a.  the Milnor fillable contact structure) and it was recently shown in \cite{cnp} that any Milnor open book on this link supports $\xi_{can}$ in the sense of Giroux \cite{gi}. Moreover, by the work of  Bogomolov and de Oliveira  \cite{bd}, the canonical contact structure of a surface singularity is Stein fillable.

\begin{figure}[ht]
  \relabelbox \small {\epsfxsize=2.5in
  \centerline{\epsfbox{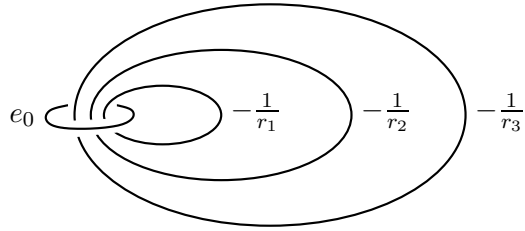}}}

\relabel{p}{$-\frac{1}{r_1}$}

\relabel{q}{$-\frac{1}{r_2}$}

\relabel{r}{$-\frac{1}{r_3}$}

\relabel{s}{$e_0$}

\endrelabelbox
        \caption{Rational surgery diagram for the small Seifert fibred $3$-manifold $Y(e_0  ; r_1, r_2, r_3)$, where $e_0 \in \bfz$ and $ r_i \in (0,1) \cap \bfq$, for
$i=1,2,3$. } \label{seif}
 \end{figure}

Let $Y = Y(e_0 ; r_1, r_2, r_3)$  denote a small Seifert fibred $3$-manifold whose
rational surgery diagram is depicted in Figure~\ref{seif}. Suppose that $Y$ is diffeomorphic to the link of some rational complex surface
singularity.  It follows that  $Y$ is a Milnor fillable rational homology sphere which is
an $L$-space (i.e., $\rk \widehat{HF}(Y) = |H_1(Y; \bfz)|$)  by a theorem of N\'{e}methi \cite{nem}.
Moreover, all tight contact structures on $Y$ are Stein fillable and they are classified, up to isotopy, in \cite{wu} (when $e_0 \leq -3$) and in \cite{gh} (when $e_0 =-2$). One of our goals in this paper is to identify
$\xi_{can}$ on $Y$,  up to isomorphism, via its Legendrian surgery diagram.

In \cite{eto2}, three numerical invariants of contact
structures were defined in terms of open books supporting the contact structures. These
invariants are the support genus $\sg(\xi)$ (the minimal genus of a
page of a supporting open book for $\xi$), the
binding number $\bn(\xi)$ (the minimal number of binding components of
a supporting open book for $\xi$ with minimal genus
pages) and the support norm  $\sn(\xi)$ (minus the maximal Euler
characteristic of a page of a supporting open book for $\xi$).

By restricting the
domain of open books to only Milnor open books,  one can redefine the above invariants specifically for the canonical contact structure $\xi_{can}$   on the link of a complex surface singularity (cf.\  \cite{ab}). We will call these invariants the Milnor
genus $\Mg(\xi_{can})$, the Milnor binding number $\Mb(\xi_{can})$, and
the  Milnor norm $\Mn(\xi_{can})$ of the canonical contact structure $\xi_{can}$.
Note that here we rather chose a simpler notation than those that were used in \cite{ab, nemet, nemtos},
to denote these invariants.  We took the liberty to introduce the notions of Milnor
genus, Milnor binding number, and Milnor norm as they have not yet appeared elsewhere in the literature.
\textit{Milnor number}, however, is well-known and it corresponds
to the first Betti number of the page of the Milnor open book, in our context.

It is clear by definition that  $\sg(\xi_{can})  \leq  \Mg(\xi_{can})$ and  $\sn(\xi_{can})  \leq  \Mn(\xi_{can})$, but no such inequality exists for the binding
numbers in general. In Section~\ref{mvs}, we describe  an infinite  family of  Milnor fillable contact
$3$-manifolds so
that $\sg(\xi_{can}) < \Mg(\xi_{can})$ and $\sn(\xi_{can}) < \Mn(\xi_{can})$. In fact, we show the existence of a sequence of Milnor fillable contact $3$-manifolds such that $\Mg(\xi_{can})  - \sg(\xi_{can})$ and
$\Mn(\xi_{can})  - \sn(\xi_{can})$ are arbitrarily large.
As a consequence,
we deduce that Milnor open books are neither norm nor genus minimizing.
We find this result interesting since there are other instances in geometric-topology, where  the ``complex representatives'' are minimizers. Most notably, the link of a
complex \textit{plane} curve singularity bounds a smooth complex curve of genus equal to its
Seifert  genus.

It is well-known \cite{or} that the small Seifert fibred $3$-manifold $Y = Y(e_0  ; r_1, r_2, r_3) $ can also be described by a star-shaped weighted plumbing tree, where the central vertex has weight $e_0$ and there are three legs emanating from that vertex corresponding to the continued fraction expansions of the rational numbers $-\frac{1}{r_i}$, for $i=1,2,3$, respectively. If $e_0 \leq -3$, then we can easily identify $\xi_{can}$ on $Y$, based on the work in \cite{eo, oz, wu} (see Section~\ref{pla}).

To handle the case $e_0=-2$, the construction in \cite{eto} can be used together with the results in \cite{b} to find explicit Milnor open books and the corresponding  $\xi_{can}$, but we ultimately identify the canonical contact structure based on the adjunction formula (see Theorem~\ref{adj}).

In this paper, we will denote a right-handed (resp. left-handed) Dehn twist along a curve $\a$ as $\a$ (resp. $\a^{-1}$), for the sake of simplicity. We refer the reader to \cite {et2} and \cite {ozst} for more on open books and contact structures and to \cite{neme} for normal surface singularities.

\section{Open books and contact structures} \label{sec:openbook}

Suppose that for an oriented link $L$ in a closed and oriented
$3$-manifold $Y$ the complement $Y\setminus L$ fibers over the
circle as $\pi \colon Y \setminus L \to S^1$ such that
$\pi^{-1}(\theta) = \Sigma_\theta $ is the interior of a compact
surface with $\partial\Sigma_\theta = L$, for all $\theta \in S^1$. Then $(L, \pi)$ is
called an \emph{open book decomposition} (or just an \emph{open
book}) of $Y$. For each $\theta \in S^1$, the surface
$\Sigma_\theta$ is called a \emph{page}, while $L$ the
\emph{binding} of the open book. The monodromy of the fibration
$\pi$ is defined as the diffeomorphism of a fixed page which is
given by the first return map of a flow that is transverse to the
pages and meridional near the binding. The isotopy class of this
diffeomorphism is independent of the chosen flow and we will refer
to that as the \emph{monodromy} of the open book decomposition.

An open book $(L, \pi)$ on a $3$-manifold $Y$ is said to be
\emph{isomorphic} to an open book $(L^\prime, \pi^\prime)$ on a
$3$-manifold $Y^\prime$, if there is a diffeomorphism $f: (Y,L) \to
(Y^\prime, L^\prime)$ such that $\pi^\prime \circ f = \pi$ on $Y
\setminus L$. In other words, an isomorphism of open books takes
binding to binding and pages to pages.

An open book can also be described as follows. First consider the
mapping torus $$\Sigma_\phi= [0,1]\times \Sigma/(1,x)\sim (0,
\phi(x))$$ where $\Sigma$ is a compact oriented surface with $r$
boundary components and $\phi$ is an element of
 the mapping class group $\Gamma_\Sigma$ of $\Sigma$.
 Since
$\phi$ is the identity map on $\partial \Sigma$,
the boundary $\partial \Sigma_\phi$ of the
mapping torus $\Sigma_\phi$
can be canonically identified with $r$ copies of $T^2 =
S^1 \times S^1$, where the first $S^1$ factor is identified with $[0,1] /
(0\sim 1)$ and the second one comes from a component of $\partial \Sigma$.
Now we glue
in $r$ copies of $D^2\times S^1$ to cap off $\Sigma_\phi$
so that $\partial D^2$ is
identified with $S^1 = [0,1] /
(0\sim 1)$ and
the $S^1$ factor in $D^2 \times S^1$
is identified with a boundary component of
$\partial \Sigma$. Thus we get a
closed $3$-manifold $Y= \Sigma_\phi \cup_{r} D^2 \times S^1 $ equipped with an open book
decomposition whose binding is the union of the
core circles $D^2 \times S^1$'s
that we glue to $\Sigma_\phi$
to obtain $Y$.
In conclusion, an element $\phi \in \Gamma_\Sigma$ determines a
$3$-manifold together with an ``abstract" open book decomposition on it.
Notice that by conjugating the monodromy $\phi$ of an open book on a 3-manifold
$Y$ by an element in $\Gamma_\Sigma$ we get
an isomorphic open book on a $3$-manifold
$Y^\prime$ which is diffeomorphic to $Y$.

Recall that a (positive) contact structure $\xi$ on an oriented $3$-manifold is the kernel of a $1$-form $\alpha$ such that
$\alpha \wedge d\alpha>0$. In this paper we assume that $\xi$ is coorientable, i.e., $\alpha$ is a
global 1--form.

It has been known  for a long time that every closed and oriented
$3$-manifold admits an open book decomposition. A new
interest in open books on $3$-manifolds arose recently from their
connection to contact structures, which we will describe very
briefly.

{\Def \label{compatible} An open book decomposition $(L,\pi)$ of a
$3$-manifold $Y$ is said to \emph{support} a contact structure $\xi$ on $Y$ if $\xi$ can be represented by a contact form
$\alpha$ such that $\alpha$ evaluates positively on the (oriented) binding $L$ and $d \alpha$ is a
symplectic form on every (oriented) page.} \\

In \cite {tw}, Thurston and Winkelnkemper associated a contact structure to  every open book.
It turns out that the contact structure they constructed is
in fact supported by the underlying open book. (Definition~\ref{compatible} was not available at the time.)
To state the converse we need a little digression.

Suppose that an open book decomposition with page $\Sigma$ is
specified by $\phi \in \Gamma_\Sigma$. Attach a $1$-handle to the
surface $\Sigma$ connecting two points on $\partial \Sigma$ to
obtain a new surface $\Sigma^{\prime}$. Let $\gamma$ be a closed
curve in $\Sigma^{\prime}$ going over the new $1$-handle exactly
once. Define a new open book decomposition with $ \phi^\prime= \phi
\circ t_\gamma \in \Gamma_{\Sigma^{\prime}} $, where $t_\gamma$
denotes the right-handed Dehn twist along $\gamma$. The resulting
open book decomposition is called a \emph{positive stabilization} of
the one defined by $\phi$. Notice that although the resulting monodromy depends on the chosen
curve $\gamma$, the $3$-manifold specified by $(\Sigma^\prime,
\phi^\prime)$ is diffeomorphic to the $3$-manifold specified by
$(\Sigma, \phi)$. A converse to the Thurston-Winkelnkemper result is given by

{\Thm [Giroux \cite {gi}] \label{giroux} Every contact $3$-manifold
is supported by an open book. Two open books supporting
the same contact structure admit a common positive stabilization.
Moreover two contact structures
supported by the same open book are isotopic.} \\

Next, we would like to briefly recall a
few different types of contact structures. An embedded disk $D\subset (Y, \xi )$ is an \emph{overtwisted disk} in the
contact $3$-manifold $(Y, \xi )$ if $\partial D=L$ is a Legendrian knot with $\tb_D(L)=0$, i.e., if the contact
framing of $L$ coincides with the framing given by the disk $D$. The contact $3$-manifold $(Y, \xi)$ is called \textit{overtwisted}
if it contains an overtwisted disk; $(Y, \xi)$ is called \textit{tight} otherwise.

A complex manifold $X$ is  called a \emph{Stein} manifold if it admits a proper
biholomorphic embedding into some $\bfc^N$. A compact $4$-manifold $W$ with nonempty
boundary $\partial W =Y$
is called a \emph{Stein domain} if there
is a Stein surface $X$ with plurisubharmonic function $\varphi\colon X\to [0, \infty ) $
such that $W= \varphi^{-1} \big([0 , t]\big)$ for some regular value $t$. So a compact
manifold with boundary (and a complex structure $J$ on its interior) is a Stein domain if it
admits a proper plurisubharmonic function $\varphi$ which is constant on the boundary.
Then the complex line distribution induced by $J$ is a contact
structure $\xi$ on $Y$. In this case we say that the contact
$3$-manifold $(Y, \xi )$ is \emph{Stein fillable} and $(W,
J)$ is a called a \emph{Stein filling} of $(Y,\xi )$. It is a deep result of
Eliashberg and Gromov \cite{eg} that every Stein fillable contact structure is tight.

\section{Legendrian surgery diagrams} \label{legknots}
Recall that a knot in a contact $3$-manifold is called Legendrian if it is everywhere tangent to
the contact planes. In order to have a better understanding of the topological constructions in the
later sections, we discuss a standard way to visualize
Legendrian knots in $S^3$ (actually  in $\bfr^3$)
equipped with the standard contact structure $\xi_{st}=\ker (dz+x\, dy)$. Consider a Legendrian knot
$L\subset (\bfr^3, \xi_{st})$ and take its \emph{front} projection, i.e., its
projection to the $yz$-plane. Notice that the projection has no vertical tangencies (since
$-\frac{dz}{dy}=x \neq \infty $), and for the same reason at a crossing the strand with
smaller slope is in front. It turns out that
$L$ can be $C^2$-approximated by a Legendrian knot for which the projection has only
transverse double points and \textit{cusp} singularities (see \cite{geig}, for example).
Conversely, a knot projection with these properties
gives rise to a unique Legendrian knot in $(\bfr^3,
\xi_{st})$ --- define $x$ from the projection as $x=-\frac{dz}{dy}$.
Since any projection can be isotoped to satisfy the above properties,
one can easily see that every knot can be isotoped to Legendrian
position.

The contact framing $\tb (L)$ of a knot $L$ can be computed as follows. Recall that we
measure the contact framing with respect to the Seifert framing in $S^3$. Define $w(L)$
(the \emph{writhe} of $L$) as the sum of signs of the double points
--- for this to make sense we need to fix an orientation on the knot,
but the answer will be independent of this choice. If $c(L)$ is the number of cusps, then the Thurston--Bennequin framing $\tb (L)$ given by
the contact structure is equal to $w(L)-\frac{1}{2}c(L)$ with respect to the framing
given by a Seifert surface.

Another invariant, the \emph{rotation number} $\rot(L)$ can be defined
by trivializing $\xi_{st}$ along $L$ and then taking the winding number of $TL$. For this
invariant to make sense we need to orient $L$, and the result will change sign when
reversing orientation. Since $H^2 (S^3; \bfz)=0$, this number will be independent of the
chosen trivialization. For the rotation number we have $\rot(L)= \frac{1}{2}(c_{d}(L)-c_{u}(L))$ where $c_d(L)$ $($resp. $c_u(L))$
denotes the number
of down $($resp. up$)$ cusps in the projection.

To describe the Stein fillable contact structures that we deal with in this paper we use Legendrian surgery
diagrams as follows: Consider the standard Stein $4$-ball $B^4$ with the induced standard contact structure on its boundary. Then attach Weinstein $2$-handles (cf.\ \cite{we})  along an arbitrary  Legendrian link in $\bdy B^4 =S^3$ to this ball. By the work of Eliashberg \cite{eli} the Stein structure on $B^4$ extends over the $2$-handles as long as
the attaching framing of each $2$-handle is given by $\tb -1$. The resulting Stein domain has an induced contact structure on its boundary which can be represented by the front projection of the Legendrian link along which we attach the $2$-handles. Such a front projection is called a \textit{Legendrian surgery diagram} (see \cite{go} for a thorough discussion). Notice that Legendrian surgery is equivalent to performing contact $(-1)$-surgery along
the given Legendrian link in the standard contact $S^3$ (see \cite{dg}, for example). To describe all Stein fillable contact structures one needs $1$-handles as well but those will not appear in our discussion.

\section{Milnor open books and the canonical contact structures}

Let $(X,x)$ be an isolated normal complex surface singularity. Fix a
local embedding of $(X,x)$ in $(\bfc^N, 0)$. Then a small sphere
$S^{2N-1}_{\epsilon} \subset \bfc^N$ centered at the origin intersects $X$
transversely,
and the complex hyperplane distribution $\xi_{can}$ on $M=X\cap
S^{2N-1}_{\epsilon}$ induced by the complex structure on $X$ is called the \emph{canonical}  contact structure.
It is known that, for sufficiently small
radius $\epsilon$, the contact manifold  is independent of $\epsilon$ and
the embedding, up to isomorphism. The $3$-manifold $M$ is called
the link of the singularity and $(M, \xi_{can})$ is called the \emph{contact boundary}
of $(X,x)$.

{\Def A contact manifold $(Y, \xi)$ is said to be
\emph{Milnor fillable} and the germ $(X, x)$ is called a
\emph{Milnor filling} of $(Y, \xi)$ if $(Y,
\xi)$ is isomorphic  to the contact boundary $(M, \xi_{can})$
of some isolated complex  surface singularity $(X, x)$.} \\

In addition, we say that a closed and oriented $3$-manifold
$Y$ is Milnor fillable if it carries a contact structure
$\xi$ so that $(Y,\xi)$ is Milnor fillable. The contact structure $\xi$ is called
the Milnor fillable contact structure.

By a theorem of Mumford \cite{mum}, if a contact 3-manifold is Milnor fillable, then it can be
obtained by plumbing oriented circle bundles over surfaces
according to a weighted graph with negative definite intersection
matrix. Conversely,
it follows from a well-known theorem of Grauert \cite{gr} that any 3-manifold that is given
by a plumbing oriented circle bundles over surfaces
according to a weighted graph with negative definite intersection
matrix is Milnor fillable. As for the uniqueness of Milnor fillable contact structures, we have the following fundamental result.

{\Thm \cite{cnp} Any closed and oriented $3$-manifold has at most one
Milnor fillable contact structure up to \textit{isomorphism}. } \\

In other words, every  Milnor fillable $3$-manifold $Y$ has a canonical contact structure $\xi_{can}$ which
is defined only up to isomorphism. Therefore, in this paper, we identify $\xi$ and $\xi_{can}$,
write $(Y, \xi_{can})$ instead of $(Y, \xi)$ and call $\xi_{can}$ as the canonical contact structure.

Since the ground-breaking result of Giroux \cite{gi}, the geometry of contact
structures is often studied
via their topological counterparts, namely the open book decompositions.  In the
realm of
surface singularities this fits nicely with some work of Milnor \cite{m}.

{\Def Given an analytic function $f\colon (X,x) \rightarrow (\bfc , 0)$
vanishing at $x$, with an isolated singularity at $x$, the open book decomposition
$\OB_f$ of the boundary $M$ of $(X,x)$ with binding $L=M \cap
f^{-1} (0)$ and projection $ \pi=\frac{f}{|f|}\colon M \setminus L \to
S^1 \subset \bfc$ is called the \emph{Milnor open book} induced by
$f$. } \\

Such functions $f$ exist  and one can talk about many Milnor open books on the singularity
link $M$. Therefore, there are many Milnor open books on any given Milnor fillable contact $3$-manifold
$(Y, \xi)$, since, by definition, it is isomorphic to the link $(M, \xi_{can})$
of some isolated complex  surface singularity $(X, x)$.

Milnor open books have two essential features \cite{cnp}: They all support
the canonical contact structure $\xi_{can}$ and they are
horizontal when considered on the plumbing description of a Milnor
fillable $3$-manifold. Recall that an open book on a plumbing of circle bundles is called
\emph{horizontal} if its binding is a collection of some fibres and
its pages are transverse to the fibres. One usually requires that the
orientation induced on the binding by the pages coincides with the
orientation of the fibres induced by the fibration.

\section{Rational surface singularities}
Let $(X,x)$ be a germ of a
normal complex surface having a singularity at $x$. Denote by $X$ a
sufficiently small representative of $(X,x)$.  Fix a resolution
$\pi\colon \tilde X\to X$ of $(X,x)$ and denote the irreducible
components of the exceptional divisor $E=\pi^{-1}(x)$ by
$\bigcup_{i=1}^r E_i$. The {\em fundamental cycle} of $E$ is
by definition the componentwise smallest nonzero effective divisor $Z=\sum z_i E_i$
satisfying $Z\cdot E_i\leq 0$ for all $i$.

{\Def[\cite{a}]  The singularity at $x$ of the germ $(X,x)$
is called {\em rational} if each irreducible component $E_i$ of the
exceptional divisor $E$ is isomorphic to $\bfc P^1$ and
$$ Z\cdot Z+\sum_{i=1}^n z_i(-E_i^2-2)= -2, $$
where $Z=\sum z_iE_i$ is the fundamental cycle of $E$.
}\\

Suppose now that $(X,x)$ is a germ of a normal complex  surface having a rational singularity at $x$. In \cite{ab}, the following result is proved:

{\Thm \label{Milmin}
Both the page-genus and the page-genus plus the number of binding components of the Milnor open book $\OB_f$ are minimized when $f$ is taken to be the restriction of a generic linear form on $\bfc^N$ to $(X,x)$ for some / any local embedding of $(X,x)$ in $(\bfc^N,0)$.
}\\

Let $\OB_{min}$ denote the Milnor open book given by taking the restriction of a generic linear form on $\bfc^N$ to $(X,x)$ for some local embedding of $(X,x)$ in $(\bfc^N,0)$. We will call $\OB_{min}$ the {\em minimal} Milnor open book. Then it is clear from Theorem \ref{Milmin} that $\Mg(\xi_{can})=g(\OB_{min})$ and $\Mb(\xi_{can})=\bc(\OB_{min})$, where
$g(\OB)$ (resp.\ $\bc(\OB)$) denotes the page-genus (resp.\ the number of binding components) of the open book $\OB$. For the Milnor norm, note that, from the definition,
$$
\Mn(\xi_{can}) = \min\{ 2g({\OB})-2+\bc({\OB}) \},
$$
where the minimum is taken over all compatible Milnor open books ${\OB}$. Hence it also follows from Theorem \ref{Milmin} that
$$ \Mn(\xi_{can}) = 2g({\OB_{min}})-2+\bc({\OB_{min}}) = 2\Mg(\xi_{can})-2+\Mb(\xi_{can}). $$

Suppose that $\pi\colon \tilde X\to X$ is a good resolution of $(X,x)$ and let $E_1,\ldots,E_r$ denote the irreducible components of the exceptional divisor $E$. Given an analytic function $f\colon (X,x) \rightarrow (\bfc , 0)$
vanishing at $x$, with an isolated singularity at $x$, the open book decomposition
$\OB_f$ is a horizontal open book with binding a vertical link of type $\underline n=(n_1,\ldots,n_r)$, where the $n_i$ are defined as follows: Consider the decomposition $(f\circ\pi)=(f\circ\pi)_e+(f\circ\pi)_s$ of the divisor $(f\circ\pi)\in\Div(\tilde X)$ into its exceptional and strict parts
given by $(f\circ\pi)_e$ is supported on $E$ and $\dim(|(f\circ\pi)_s\cap E|<1$. It can be shown that the $r$-tuple $\underline n=(n_1,\ldots,n_r)$ satisfies
\begin{equation} \label{incid}
 I(\G)\underline m^t=-\underline n^t
\end{equation}
for some $r$-tuple $\underline m=(m_1,\ldots,m_r)$ of positive integers, where $\G$ denotes the intersection matrix of $E$.

On the other hand, it follows from Artin \cite{a} that for any $r$-tuple $\underline n$ of nonnegative integers which satisfies \eqref{incid} for some $r$-tuple $\underline m$ of positive integers there is a Milnor open book decomposition of the boundary of $(X,x)$ whose binding is equivalent to a vertical link of type $\underline n$. It can be shown that if $Z=\sum_{i=1}^r z_i E_i$ is the fundamental cycle of the resolution $\pi$, then the above construction for the $r$-tuple $\underline m=(z_1,\ldots,z_r)$ gives the minimal Milnor open book $\OB_{min}$.

{\Rem In \cite{nemtos}, a generalization of Theorem \ref{Milmin} is given for all Milnor fillable rational homology $3$-spheres.}

\section{Planar Milnor open books} \label{pla}

A vertex in a weighted plumbing graph is called a bad vertex if the sum of the Euler number $e$ and the degree $d$  of that vertex is positive.

{\Prop \label{nobad}  Let $Y$ be the link of a rational surface singularity presented by a weighted plumbing tree  with $r$ vertices.  Assume that there are no bad vertices in the plumbing tree. Then we have $\Mg(\xi_{can})= 0$,  $\Mb(\xi_{can}) = -\displaystyle\sum_{i=1}^r   (e_i+d_i)$, and $\Mn(\xi_{can})= -2 + \Mb(\xi_{can})$.  }

\begin{proof}

The link $Y$ admits a planar horizontal open book $\OB$ (cf.\ \cite{eo}) with binding a vertical link of type $$\underline
n=-(e_1+d_1, e_2+d_2, \ldots, e_r+d_r)$$ so that for  $\underline m=(1,1, \ldots, 1)$  we have
$$I(\G)\underline m^t=-\underline n^t,$$ where $I(\G)$ denotes the intersection matrix of the tree $\G$ which defines $Y$. Suppose that $Y$ is the link of the rational surface singularity $(X,x)$. Then $\underline m$ corresponds to the fundamental cycle of the minimal resolution of $(X,x)$. Since the open book $\OB$ and the minimal Milnor open book $\OB_{min}$ on the rational homology $3$-sphere $Y$ have equivalent bindings, it follows,
by a result of Caubel and Popescu-Pampu \cite{cp}, that $\OB$ is isotopic to $\OB_{min}$. This proves that $\Mg(\xi_{can})=g(\OB)=0$, $\Mb(\xi_{can})=\bc(\OB)=-\displaystyle\sum_{i=1}^r   (e_i+d_i)$ and
$\Mn(\xi_{can})=2\Mg(\xi_{can})-2+\Mb(\xi_{can})=-2+\Mb(\xi_{can})$.
\end{proof}

{\Rem The canonical
contact structure $\xi_{can}$ of the link of a singularity as in Proposition~\ref{nobad} can be explicitly given
by a Legendrian surgery diagram using the methods discussed in \cite{oz}. } \\

The links of such singularities include lens spaces.  Recall that the lens space $L(p,q)$ is obtained from $S^3$ by $-p/q$ surgery on the unknot.  Let $ [ a_1, a_2, \ldots, a_n  ]$ denote the continued fraction expansion of the rational number $-p/q$, where  $a_i \leq -2$, for  $ i=1,2, \ldots, n$.  The next result immediately follows from Proposition~\ref{nobad}. (Note that a Legendrian surgery diagram for  $\xi_{can} $ on $L(p, q)$ is given
by Figure 4 in \cite{oz}.)

{\Cor  For the canonical contact structure $\xi_{can} $ on $L(p, q)$, we have  $\Mg(\xi_{can})= 0$,  $\Mb(\xi_{can}) = 2-2n -\displaystyle\sum_{i=1}^n  a_i$, and $\Mn(\xi_{can})= -2 + \Mb(\xi_{can})$.} \\

In particular,   the Milnor binding number and hence the Milnor norm can be made arbitrarily large by choosing, say, $a_1$ arbitrarily
small, for fixed $n$. It is known (cf. \cite{sch}) that the support genus is zero for all (tight) contact structures on  lens spaces.  We would like to conjecture that
 $$ \bn (\xi_{can}) = \Mb(\xi_{can})  \;  \mbox{and} \;  \sn (\xi_{can}) = \Mn(\xi_{can}), $$  for the canonical contact structure $\xi_{can} $ on $L(p, q)$. This is certainly true for $L(n, n-1)$, for $n \geq 2$. \\

The link of singularities in Proposition~\ref{nobad} also include small Seifert fibred manifolds of the form $Y=Y(e_0; r_1,r_2,r_3)$ where $e_0 \leq -3$. Note that in the star-shaped plumbing tree of $Y$ all the weights are less than or equal to $-2$ and therefore one can Legendrian realize these unknots to obtain distinct Stein fillable contact
structures on $Y$. The construction in \cite{oz}, coupled with Wu's classification \cite{wu}, allows us to conclude that

{\Prop  The canonical contact structure $\xi_{can}$ on $Y=Y(e_0; r_1,r_2,r_3)$, where $e_0 \leq -3$, can be identified as the contact structure obtained by putting all the extra zigzags of the Legendrian unknots in the star-shaped presentation of $Y$, to one fixed side.} \\

Moreover, an explicit planar Milnor open book for $\xi_{can}$ can be
described as in the proof of Proposition~\ref{nobad}.

\section{Elliptic Milnor open books}

In this section, we examine some families of rational surface singularities whose links are diffeomorphic to
some small Seifert fibred $3$-manifold of the form $Y(-2; r_1,r_2,r_3)$.

\subsection{First family}

\begin{figure}[ht]
  \relabelbox \small {\epsfxsize=5in
  \centerline{\epsfbox{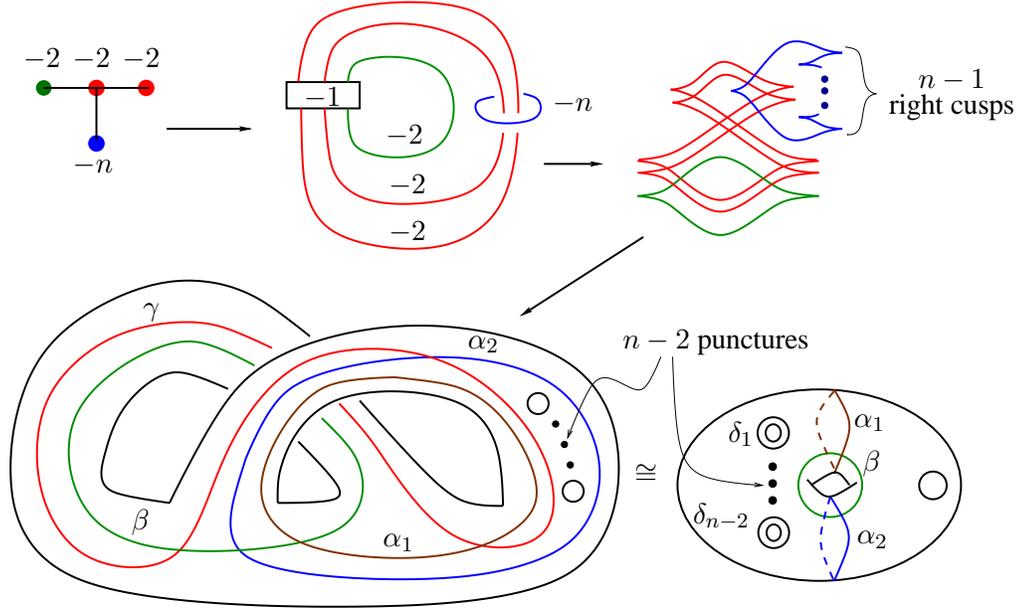}}}

\relabel{1}{$\a_1$}

\relabel{5}{$\a_2$}

\relabel{3}{$\b$}

\relabel{4}{$\g$}

\relabel{6}{$-1$}

\relabel{7}{$-2$}

\relabel{8}{$-2$} \relabel{cong}{$\cong$}

\relabel{9}{$-2$}

\relabel{10}{$-n$}

\relabel{11}{$\a_1$}

\relabel{12}{$\a_2$}

\relabel{13}{$\b$}

\relabel{a}{$-2$} \relabel{n}{$n-2$ punctures}

\relabel{m}{$\;\;\;\;n-1$ }
\relabel{o}{right cusps}

\relabel{b}{$-2$}

\relabel{c}{$-2$}

\relabel{d}{$-n$}
\relabel{e}{$\d_1$} \relabel{f}{$\d_{n-2}$}

\endrelabelbox
        \caption{From plumbing tree to a particular open book on $N_n$} \label{plumbing}
 \end{figure}

Let  $N_{n} = Y(-2; \frac{1}{2},\frac{1}{2},\frac{1}{n})$, for $n \geq 2$,  which is given by the plumbing diagram in Figure~\ref{plumbing}. Note that this resolution graph comes from a rational surface singularity. Now we apply the recipe described in \cite{eto} to obtain a Legendrian surgery diagram for a particular Stein fillable contact structure on $N_n$. First we roll up the plumbing diagram  and then Legendrian realize the given surgery curves, by paying attention to their Thurston-Bennequin invariants as shown in Figure~\ref{plumbing}.

Next, we construct an open book of $N_n$ supporting this particular contact structure: We start from an open book of $S^3$ and then embed the surgery curves onto the pages of this open book. The initial page is a torus with one boundary component and the monodromy of this open book of $S^3$ before the surgery is given by  $ \b \a_1$, where $\a_1$ and $\b$ generate the first homology group of the page. Now we apply the Legendrian surgeries along the four given curves
to get an open book of $N_n$ with the  monodromy $\phi_n= \a_2 \g^2 \b^2 \a_1 \d_1 \d_2 \cdots \d_{n-2}$, where $\d_1, \d_2, \ldots, \d_{n-2}$ are boundary parallel curves around the $n-2$ small punctures on the torus, respectively. These punctures occur as a result of
stabilizing the page appropriately. Next we apply some mapping class group tricks. Move $\b$
over $\g^2$ to the left and use the fact that $\g\b = \b \a_1$ to get  $ \phi_n = \a_2 \b \a_1 \a_1 \b \a_1  \d_1 \d_2 \cdots \d_{n-2} $. Then use the well-known the braid relations  and some simple overall conjugations
to obtain
\begin{align*}  \phi_n &= \a_2 \b \a_1 \b \a_1 \b  \d_1 \d_2 \cdots \d_{n-2} \\
&= \b \a_2 \b \a_1 \b \a_1  \d_1 \d_2 \cdots \d_{n-2} \\
&= \a_2 \b \a_2 \a_1 \b \a_1  \d_1 \d_2  \cdots \d_{n-2} \\
&= \a_1 \a_2 \b \a_1 \a_2 \b  \d_1 \d_2 \cdots \d_{n-2}  \\
&= (\a_1 \a_2 \b)^2 \d_1 \d_2\cdots \d_{n-2}.
\end{align*}

\begin{figure}
\begin{center}
\includegraphics[scale=1.0]{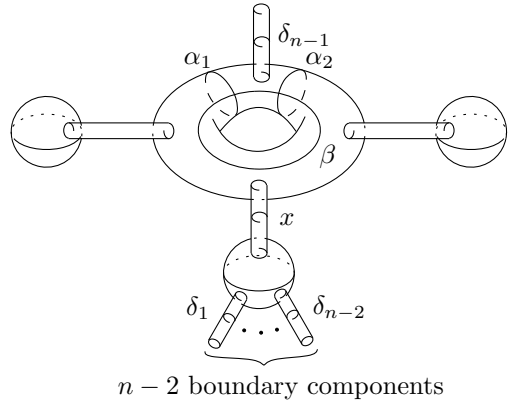}
\end{center}
\caption{A page of a Milnor open book on $N_n$}
\label{milnorpage}
\end{figure}

We claim that the open book we just constructed is isomorphic to a Milnor open book. To prove our claim, we construct a Milnor open book on $N_n$ using the method of \cite{b}.
(Note that the corresponding singularity is a dihedral singularity and
an explicit presentation of an associated Milnor open book is also given in \cite{y}.) The construction is as follows:
Enumerate the
vertices of the plumbing graph from left to right along the top row with the $-n$ vertex coming last
and consider the 4-tuple of positive integers $\underline
m=(1,2,1,1)$. Note that this corresponds to the fundamental cycle of the minimal resolution of the associated rational singularity and hence the Milnor open book that it corresponds to is the minimal Milnor open book $\OB_{min}$. The page of $\OB_{min}$ is depicted in Figure~\ref{milnorpage}.

From the construction in \cite{b}, the monodromy $\psi_n$ of $\OB_{min}$ satisfies
$$ \psi_n^2=\d_{n-1} x \d_1^2 \d_2^2 \cdots \d_{n-2}^2. $$
We now use the fact that elements in mapping class groups of surfaces of genus one with nonempty boundary can have at most one $m$-root up to conjugation (cf.\ \cite{bp}) to deduce that $\psi_n$ may be written
$$ \psi_n=(\a_1 \a_2 \b)^2 \d_1 \d_2 \cdots \d_{n-2}$$ which indeed coincides with $\phi_n$ above.
Here we used the two-holed torus relation $$\d_{n-1} x = (\a_1 \a_2 \b)^4.$$

\begin{figure}[ht]
  \relabelbox \small {\epsfxsize=2in
  \centerline{\epsfbox{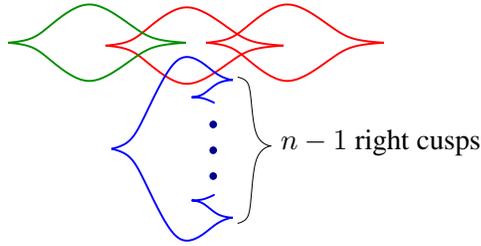}}}
\relabel{m}{$n-1$ right cusps}

\endrelabelbox
        \caption{The canonical contact structure $\xi_{can}$ on $N_n$} \label{leg}
 \end{figure}

It follows that the contact structure depicted by  the Legendrian surgery diagram in Figure~\ref{plumbing}
is isomorphic to the canonical contact structure. By the classification of the tight contact structures on $N_n$ (see \cite{gh}), this contact structure is in fact isomorphic to the contact structure whose Legendrian surgery diagram is
given as in Figure~\ref{leg}.

Note that there are $n-1$ distinct tight contact structures on $N_n$, (all of which are Stein fillable),  obtained by distributing the  zigzags of  the bottom curve in different ways.  The one where one puts all the extra zigzags to the left  is also isomorphic to $\xi_{can}$.

When $n=2$, the page is a once punctured torus, $\a_1 = \a_2 =\a$ and the monodromy of the open book
is given by  $$ \phi_2 =  (\a^2 \b)^2 = (\a\b)^3, $$
which coincides with the  Milnor open book presented in \cite{b}.

\subsection{Second family}

\begin{figure}[ht]
  \relabelbox \small {\epsfxsize=5in
  \centerline{\epsfbox{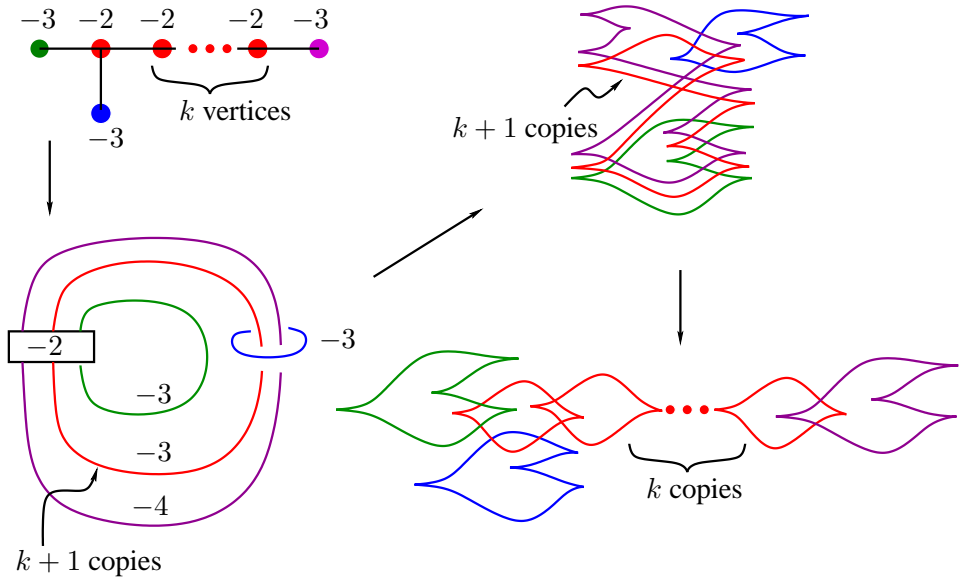}}}

\relabel{6}{$-2$}

\relabel{7}{$-3$}

\relabel{8}{$-3$}

\relabel{9}{$-4$}

\relabel{10}{$-3$}

\relabel{a}{$-3$} \relabel{b}{$-2$}

\relabel{c}{$-2$} \relabel{d}{$-3$}

\relabel{e}{$-3$} \relabel{f}{$-2$}

\relabel{p}{$k+1$ copies}

\relabel{q}{$k+1$ copies}

\relabel{r}{$k$ vertices}

\relabel{s}{$k$ copies}

\endrelabelbox
        \caption{The canonical contact structure $\xi_{can}$ on $M_k$, for
$k \geq 0 $ } \label{fam}
 \end{figure}

\begin{figure}[ht]
  \relabelbox \small {\epsfxsize=5in
  \centerline{\epsfbox{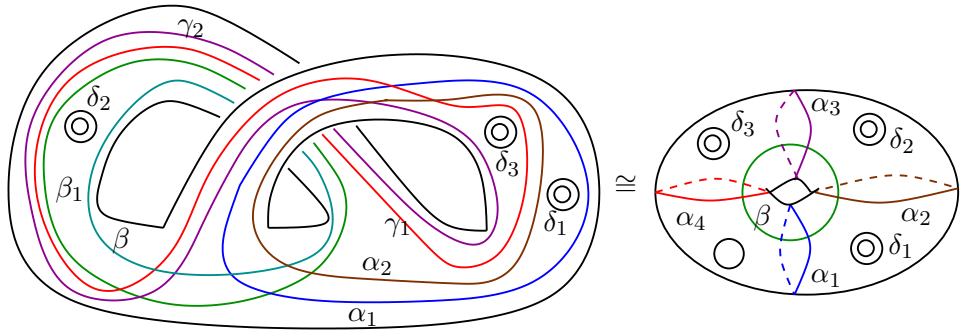}}}

\relabel{1}{$\a_2$}

\relabel{5}{$\a_1$}

\relabel{3}{$\b$}

\relabel{2}{$\g_2$}

\relabel{4}{$\g_1$}

\relabel{b2}{$\b_1$}

\relabel{11}{$\a_3$}

\relabel{12}{$\a_1$}

\relabel{13}{$\b$} \relabel{14}{$\a_2$} \relabel{15}{$\a_4$}

\relabel{cong}{$\cong$}

\relabel{a}{$\d_1$} \relabel{c}{$\d_2$} \relabel{b}{$\d_3$}

\relabel{e}{$\d_1$} \relabel{f}{$\d_2$} \relabel{g}{$\d_3$}

\endrelabelbox
        \caption{Page of a particular open book on $M_k$} \label{plumbing2}
 \end{figure}

Let $M_k=Y(-2; \frac{1}{3},\frac{1}{3},\frac{2k+1}{2k+3})$, for $k \geq 0$ as depicted in Figure~\ref{fam}.   Note that this resolution graph comes from a rational surface singularity. First we roll up the plumbing diagram and then Legendrian realize the given surgery curves in a certain way to obtain a particular Stein fillable contact structure on $M_k$. Next, we construct an open book of $M_k$ supporting  this contact structure by the method discussed in \cite{eto}: We start from an open book of $S^3$, where the page is a torus with one boundary component and the monodromy is given by
$\b \a_2$. Then embed the surgery curves onto the pages as depicted in  Figure~\ref{plumbing2}.

The page of the open book on $M_k$ is a torus with four punctures (see Figure~\ref{plumbing2})
and the monodromy of the open book after the surgery is given by
\begin{align*}  \phi_k &= \a_1 \g_2 \g_1^{k+1} \b_1 \b \a_2 \d_1 \d_2 \d_3  \\
    &= \a_1 \b  \a_4  \a_3^{k+1} \b_1 \a_2 \d_1 \d_2 \d_3 \\
      &= \a_1 \b \a_4  \a_3^{k+1}  \b_1 \a_2 \a_3^{-1} \a_3 \d_1 \d_2 \d_3 \\
        &=  \a_1 \b \a_4  \a_3^{k+1} \a_2 \a_3^{-1} \b  \a_3 \d_1 \d_2 \d_3 \\
   &=  \a_1 \b \a_4  \a_3^k \a_2  \b  \a_3 \d_1 \d_2 \d_3 \\
 &=  \a_1 \b \a_4  \a_2  \a_3^k \b  \a_3 \d_1 \d_2 \d_3.
\end{align*}

Note that we moved $\b$ over $\g_1$ and $\g_2$ to the left in the second equation above, naming
the resulting curves $\a_3$ and $\a_4$, respectively. Also we moved
$\a_2 \a_3^{-1}$ to the left over $\b_1$ in the fourth equation. By using braid relations and simple induction we obtain that
$$ \phi_k= \a_1 \a_3 \b \a_2  \a_4 \a_1^s \a_3^{k-s} \b \d_1 \d_2 \d_3\quad \text{for any } 0 \leq s \leq k. $$

Rename $\a_i$ as $\a_{i+1}$ mod $4$, to cyclically permute the $\a$ curves on the torus.
It follows that
$$\phi_k = \begin{cases}
 \d_1 \d_2 \d_3 \a_2^{(k-1)/2} \a_4^{(k-1)/2} \a_1 \a_3 \b \a_2 \a_4 \b  \a_2  & \text{when $k$ is odd (and $s=(k+1)/2$), } \\
 \d_1 \d_2 \d_3 \a_2^{k/2} \a_4^{k/2} \a_1 \a_3 \b \a_2 \a_4 \b    &  \text{when $k$ is even (and $s=k/2$).}
\end{cases} $$

\begin{figure}
\begin{center}
\includegraphics[scale=1.0]{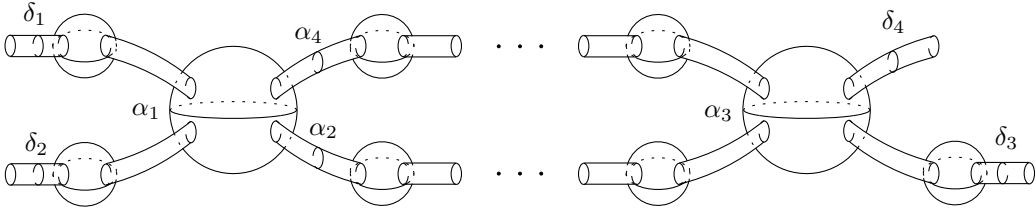}
\end{center}
\caption{A Milnor open book for $M_k$}
\label{milnorpagefam}
\end{figure}

Now we claim that this particular open book on $M_k$ agrees with a Milnor open book. To see this, enumerate the vertices of the plumbing graph from left to right along the top row with the bottom vertex coming last and consider the $(k+4)$-tuple of positive integers $\underline m=(1,2,\ldots,2,1,1)$. The page of the associated Milnor open book, which is $\OB_{min}$, is depicted in Figure~\ref{milnorpagefam}.

From the construction in \cite{b}, the monodromy $\psi_k$ satisfies
$$ \psi_k^2=\d_1^3\d_2^3\d_3^3\d_4\a_2^{k}\a_4^{k}. $$
We claim $\psi_k=\phi_k$.
To see this, using the uniqueness result from \cite{bp}, clearly it is sufficient to check that
$$ (\d_1\d_2\d_3\a_1\a_3\b\a_2\a_4\b \a_2^s\a_4^{k-s})^2=\d_1^3\d_2^3\d_3^3\d_4 \a_2^k\a_4^k $$
for any $0\leq s\leq k$. Cancelling $\d_1^2\d_2^2\d_3^2$ on both sides and using the four-holed torus relation
$(\a_1\a_3\b\a_2\a_4\b)^2=\d_1\d_2\d_3\d_4$ (cf.\ \cite{ko}),
we thus need to check
$$ \a_1\a_3\b\a_2\a_4\b \a_2^s\a_4^{k-s} \a_1\a_3\b\a_2\a_4\b \a_2^s\a_4^{k-s}= (\a_1\a_3\b\a_2\a_4\b)^2\a_2^k\a_4^k.  $$ Making further cancellations, it is thus sufficient to check that
$$ \a_2^s\a_4^{k-s} \b\a_2\a_4\b \a_2^s\a_4^{k-s}=\b\a_2\a_4\b \a_2^k\a_4^k $$
for any $0\leq s\leq k$, which can be done by using the braid relations and a simple induction argument.

Note that there are exactly eight (cf.\ \cite{gh}) distinct tight contact structures on $M_k$. The contact structure in Figure~\ref{fam} is isomorphic to $\xi_{can}$ on $M_k$. Putting all the extra zigzags to the left rather than to the right in  Figure~\ref{fam} would yield an isomorphic contact structure.



\section{Milnor versus support genus}\label{mvs}

\begin{figure}[ht]

  \relabelbox \small {\epsfxsize=5.5in
  \centerline{\epsfbox{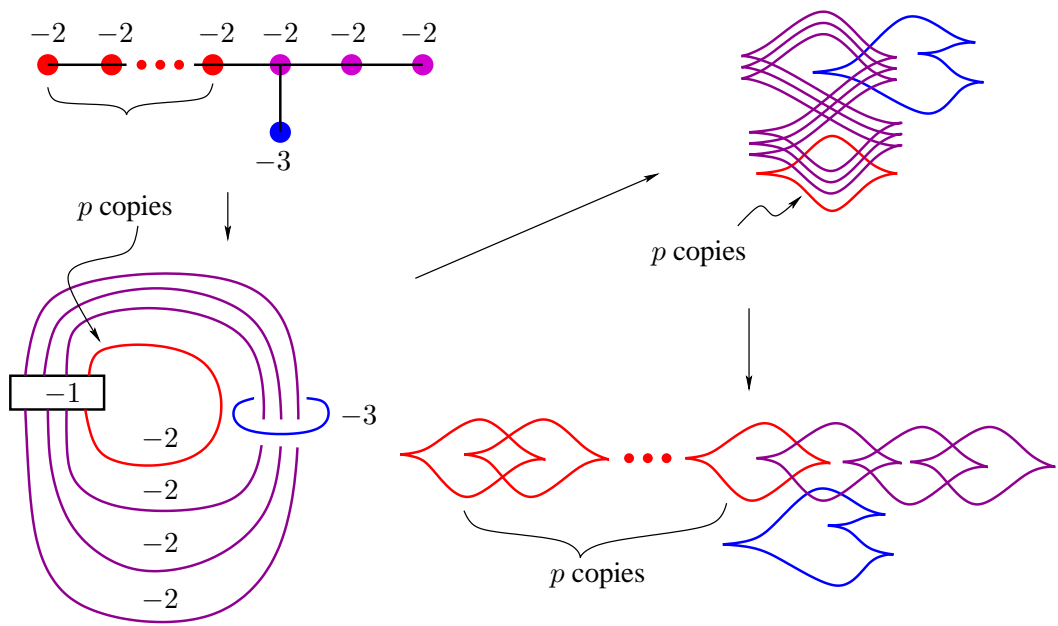}}}

\relabel{6}{$-1$}

\relabel{7}{$-2$}

\relabel{8}{$-2$}

\relabel{9}{$-2$} \relabel{11}{$-2$}

\relabel{10}{$-3$}

\relabel{a}{$-2$}

\relabel{b}{$-2$}

\relabel{c}{$-2$}

\relabel{d}{$-3$}

\relabel{e}{$-2$} \relabel{f}{$-2$} \relabel{g}{$-2$}

\relabel{p}{$p$ copies}

\relabel{q}{$p$ copies}

\relabel{s}{$p$ copies}

\endrelabelbox
           \caption{The contact structure $\xi_1 \cong \xi_{can}$ on $Y_p$} \label{yp}
 \end{figure}

In this section, we describe  an infinite  family of  Milnor fillable contact
$3$-manifolds so
that $\sg(\xi_{can}) < \Mg(\xi_{can})$ and $\sn(\xi_{can}) < \Mn(\xi_{can})$.  Consider the small Seifert fibred $3$-manifold $Y_p = Y(-2  ; \frac{1}{3}
, \frac{2}{3}, \frac{p}{p+1})$, for $p \geq 2$.   First, we observe that $Y_p$ is (diffeomorphic to) the
link of a rational complex surface singularity, and
its resolution graph $\G_p$ is shown in Figure~\ref{yp}.
By the classification of the tight contact structures on $Y_p$ given by Ghiggini
\cite{gh}, there are exactly two nonisotopic tight contact structures $\xi_1$ and
$\xi_2$ on  $Y_p$, both of which are Stein fillable.

 {\Prop For $i=1,2$, we have  $\sg(\xi_i) \leq 1$ and $\sn(\xi_i)=2$.}

 \begin{proof}

According to the recipe in
\cite{eto}, first we roll up the plumbing diagram (i.e., we appropriately slide
handles) and then Legendrian realize the given surgery
curves in a certain way to obtain the Legendrian surgery diagram for a particular
Stein fillable contact
structure $\xi_1$ on $Y_p$ (see Figure~\ref{yp}).
\begin{figure}[ht]
  \relabelbox \small {\epsfxsize=5in
  \centerline{\epsfbox{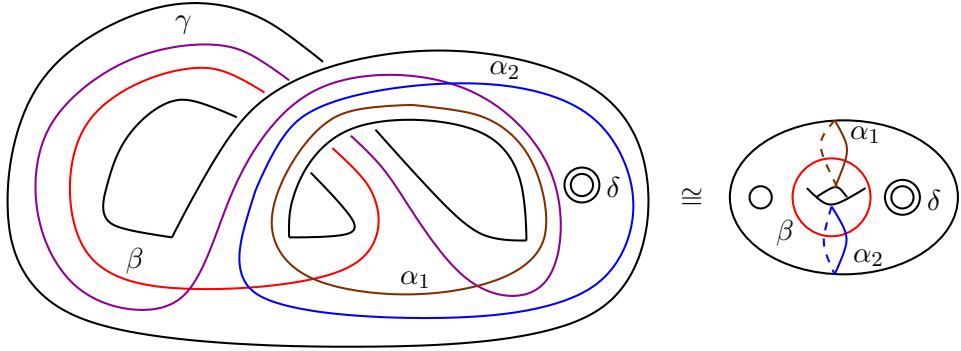}}}

\relabel{1}{$\a_1$}

\relabel{5}{$\a_2$}

\relabel{3}{$\b$}

\relabel{4}{$\g$}

\relabel{cong}{$\cong$}

\relabel{11}{$\a_1$}

\relabel{12}{$\a_2$}

\relabel{13}{$\b$}

\relabel{e}{$\d$} \relabel{f}{$\d$}

\endrelabelbox
        \caption{The page of an open book compatible with $\xi_1$} \label{openbook}
 \end{figure}

Next, we construct
an open book of $Y_p$ supporting this contact structure: We start from an open
book of $S^3$ and
then embed the surgery curves onto the pages as depicted in  Figure~\ref{openbook}.  The initial page is a torus with one
boundary component and the monodromy of this open book of $S^3$ before the surgery
is given by  $ \b \a_1$, where $\a_1$ and $\b$ generate the first homology group of
the page. Now we apply the Legendrian surgeries along the given curves
to get an open book of $Y_p$ with the  monodromy $\phi_p= \a_2 \g^3 \b^p \b \a_1 \d
$, where $\d$ is parallel to the small puncture on the torus, which occurs as a
result of
stabilizing the page appropriately. Next we apply some mapping class group tricks.
Move $\b$
over $\g^3$ to the left and use the fact that $\g\b = \b \a_1$ to get  $ \phi_p =
\a_2 \b \a_1^3 \b^p \a_1  \d $. Then use the well-known the braid relations  and
some simple overall conjugations
to obtain a more symmetrical presentation of the monodromy
$$\phi_p = (\a_2 \b)^2 (\a_1 \b)^2   \b^{p-2} \delta.$$

What we described here is an  \emph{abstract}  open book which is compatible  with
$\xi_1$---where the page is a torus with two boundary
components  and monodromy  is $\phi_p$.  Let $\xi_2$ denote the contact structure
where we put the extra zigzag
in Figure~\ref{yp} to the left. Note that one can not distinguish the abstract open
books corresponding to $\xi_1$ and $\xi_2$ and in fact
$\xi_1$ is isomorphic to  $\xi_2$. It
follows that $\sg(\xi_i) \leq 1$, since we have already constructed a genus one open
book compatible with  $\xi_i$.

One can show that
$$H_1 (Y_p, \bfz) = \begin{cases}
\bfz_3  \bigoplus  \bfz_3  & \text{$p = 2$  mod $3$, } \\
 \bfz_9  &  \text{otherwise.}
\end{cases} $$ \\

It turns out that, for all $p$ congruent to $2$ mod $3$, Poincar\'{e} dual
$PD (e (\xi_i)) \in H_1 (Y_p, \bfz)$  of the Euler class $e (\xi_i)$ is a generator of one of the
$\bfz_3$-factors. Similarly $PD (e (\xi_i))$   is a generator of $ H_1 (Y_p, \bfz)$
when $p$ is not congruent to $2$ mod $3$. Therefore the contact structure $\xi_i$
can not be compatible
with an elliptic open book  with connected binding by  Lemma 6.1 in \cite{eto2}, since $e(\xi_i) \neq 0$.
Note that $e (\xi_1 ) =  - e(\xi_2)$, which implies that $\xi_1$ is not homotopic to  $\xi_2$ as oriented plane fields, although they are isomorphic to each other. In fact, $\xi_2$ is obtained by $\xi_1$ by reversing the orientation of the underlying plane field.

Now we claim that  $\sn(\xi_i) = 2$.  To prove our claim we need to
exclude the possibility that $\xi_i$ is compatible with a planar open book with less
than four binding components.
Suppose that $\xi_i$ is compatible with a planar open book, i.e, $\sg(\xi_i)=0$. If $\bn(\xi_i) \leq 2$,
then $\xi_i$ is the unique tight contact structure on the lens space $L(n, n-1)$ for some
$n \geq 0$ (cf.\ \cite{eto2}) which is indeed impossible since $Y_p$ is not a lens space.

Next we rule out the possibility that $\bn(\xi_i)=3$. Let $\Sigma$ be the planar surface with three boundary
components. Any diffeomorphism of $\Sigma$ is determined by three
numbers $q,r,s,$ that give the number of Dehn twists on curves
$\tau_1, \tau_2, \tau_3$ parallel to each boundary component.
Let $Y_{q,r,s}$ be the 3-manifold determined by the open book with
page $\Sigma$ and monodromy given by
$\tau_1^q \tau_2^r \tau_3^s$. It is easy to see that
$Y_{q,r,s}$ is the Seifert fibred space shown in
Figure~\ref{fig:bn3}.

\begin{figure}[ht]
  \relabelbox \small {\epsfxsize=4.5in
  \centerline{\epsfbox{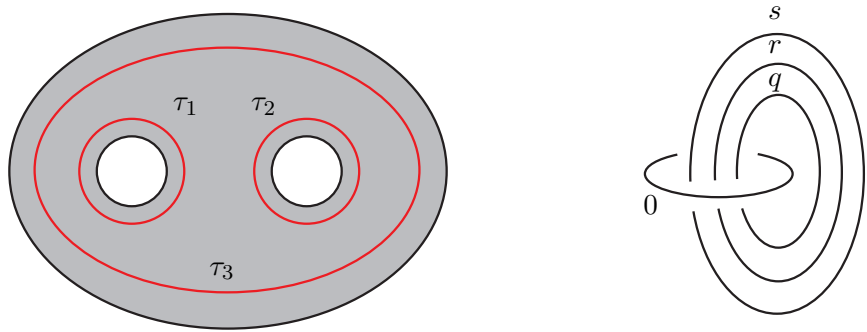}}}
  \relabel{1}{$\tau_1$}
  \relabel{2}{$\tau_2$}
  \relabel{3}{$\tau_3$}
  \relabel{0}{$0$}
  \relabel{m}{$q$}
  \relabel {n}{$r$}
  \relabel {k}{$s$}
  \endrelabelbox
        \caption{The surface $\Sigma$ on the left; the manifold $Y_{q,r,s}$ on the right.}
        \label{fig:bn3}
\end{figure}

One can compute that $ | H_1(Y_{q,r,s}, \bfz) | = qr+qs+rs$. Suppose that $\xi_i$ is compatible
with an open book with page $\Sigma$ and monodromy $\tau_1^q \tau_2^r \tau_3^s$. The tightness of
$\xi_i$ implies that the integers $q,r$, and $s$ are all nonnegative, because otherwise $\tau_1^q \tau_2^r \tau_3^s$ is not right-veering \cite{hkm}.
 Moreover, since the order of the
first homology group of $Y_p$ is $9$,  for all $p \geq 2 $,  we conclude that $(q,r,s)$ is equal to either $(0,1,9)$, $(0,3,3)$
or $(1,1,4)$.  Hence $Y_p$ is diffeomorphic to either $L(9,8)$, $L(3,2) \# L(3,2)$ or $L(9,4)$, which is a contradiction. Hence, $\bn(\xi_i) \geq 4$.  This finishes the proof of our claim that $\sn(\xi_i) = 2$.
\end{proof}

One can ask whether or not $\sg(\xi_i)=1$, although it is not essential for the purposes of this paper.
There are two known methods (cf.\ \cite{et1, oss}) of finding
obstructions to planarity of  a contact structure, but unfortunately both fail in our case.  That is because $Y_p$ is an $L$-space
and it is not an integral homology sphere.

{\Prop  We have  $\Mg(\xi_{can})=2$ and  $\Mn(\xi_{can})=3$.}

\begin{proof}
Enumerate the vertices of the plumbing graph from left to right along the top row with the bottom vertex coming last. It is then easy to check that the $(p+4)$-tuple of positive integers
$\underline m$ corresponding to the fundamental cycle of the minimal resolution of the singularity of which $Y_p$ is the link is given by $\underline m=(1,2,3,3,\ldots,3,3,2,1,1)$. The construction in \cite{b} now gives an open book decomposition $\OB({\underline m})=\OB_{min}$ of $Y_p$ with binding a vertical link of type $\underline n=(0,0,1,0,\ldots,0)$. Note that here $\underline m$ and $\underline n$ are related by
$$ I(\G_p)\underline m^t=-\underline n^t. $$

Using the formula
\begin{equation*}
g(\OB(\underline m))=1+\sum_{i=1}^{p+4} \frac{(v_i-2)m_i + (m_i-1)n_i}{2}
\end{equation*}
given in Lemma~3.1 in \cite{ab}, where $\underline m=(m_1,\ldots,m_{p+4})$ and $\underline n=(n_1,\ldots,n_{p+4})$, one has $\Mg(\xi_{can})=g(\OB(\underline m))=2$ for the unique Milnor fillable contact structure $\xi_{can}$ on each 3-manifold $Y_p$. Also one has $\Mb(\xi_{can})=\bc(\OB(\underline m))=\sum_{i=1}^{p+4}n_i=1$. It follows that $\Mn(\xi_{can})=3$, completing the proof of the proposition.

\end{proof}

Now since any Milnor
fillable contact structure is Stein fillable \cite{bd},
$\xi_{can}$ has to be isomorphic to $\xi_i$ by Ghiggini's classification \cite{gh}. Note that it does not make sense to distinguish $\xi_1$ and $\xi_2$ here since they
are isomorphic to each other. Thus

{\Cor We have $\sg(\xi_{can}) < \Mg(\xi_{can})$ and $\sn(\xi_{can}) <
\Mn(\xi_{can})$.} \\

Note, however,  that  $\Mb(\xi_{can})=1$ while  $ \bn(\xi_{can}) \geq 2$. \\

The arguments above can be generalized to prove the following

{\Thm \label{arbit} For each positive integer $k$, there is a Milnor fillable contact $3$-manifold such that $\Mg(\xi_{can})  - \sg(\xi_{can}) \geq k$ and
$\Mn(\xi_{can})  - \sn(\xi_{can}) \geq  k$.} \\

\begin{proof}
Consider the small Seifert fibred $3$-manifolds $P_n=Y(-2;\frac 1{n+1},\frac n{n+1},\frac n{n+1})$, for $n\geq 2$. Each of these is (diffeomorphic to) the link of a rational complex surface singularity with minimal resolution graph the weighted tree shown in Figure \ref{longgraph}. By the classification of tight contact structures on $P_n$ (\cite{gh}) there are exactly $n$ nonisotopic tight contact structures $\xi_i$, for $i=1,\ldots,n$, on $P_n$, each of which is Stein fillable. By \cite{eto}, for each of these tight contact structures $\xi_i$ we can find a genus one  supporting open book. Since the canonical contact structure $\xi_{can}$ on $P_n$ is tight it must be isomorphic to one of the contact structures $\xi_i$. This proves that $\sg(\xi_{can})\leq 1$. It is easy
to see that $\sn(\xi_{can}) \leq n$, since each $\xi_i$ is supported by an elliptic open book with $n$
binding components.

\begin{figure}
\begin{center}
\includegraphics[scale=1.0]{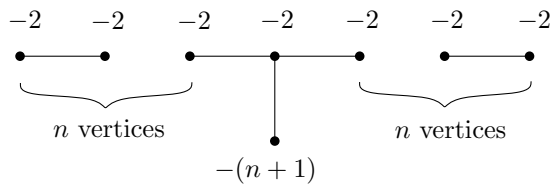}
\end{center}
\caption{The plumbing graph for $P_n$}
\label{longgraph}
\end{figure}

Now enumerate the vertices of the graph in Figure \ref{longgraph} as before and consider the $(2n+2)$-tuple
$\underline m=(1,2,3,\ldots,n-1,n,n+1,n,n-1,\ldots,3,2,1,1)$ of positive integers. This corresponds to the fundamental cycle
of the minimal resolution of the singularity of which $P_n$ is the link. It now follows that $\Mg(\xi_{can})=n$, $\Mn(\xi_{can})=2n-1$, $\sg(\xi_{can}) \leq 1$, and $\sn(\xi_{can}) \leq n$, proving the theorem, by taking $k=n-1$.
\end{proof}

{\Cor Not every open book on the link of a surface singularity is isomorphic to a Milnor open book.
}

\section{Legendrian surgery diagrams for canonical contact structures}

This section is devoted to proving the following result

{\Thm \label{adj} Let $Y = Y(e_0 ; r_1, r_2, r_3)$ be a small Seifert fibred $3$-manifold which is diffeomorphic to
the link of some rational surface singularity.
Then the canonical contact structure $\xi_{can}$
on $Y$ is given, up to isomorphism,  by the Legendrian surgery diagram
(obtained from the plumbing tree) where one puts all the extra zigzags of the Legendrian unknots to one fixed side.}

\begin{proof} Suppose that $(Y,\xi_{can})$ is diffeomorphic to the link of the rational
surface singularity $(X,x)$. Then the minimal resolution $\pi\colon \tilde X\to X$ provides a holomorphic
filling $(W,J)$ of $(Y,\xi_{can})$. In particular, $W$ is a regular neighbourhood of the exceptional
divisor $E=\bigcup E_j$ of $\pi$. Since the curves $E_j$ are holomorphic, by the adjunction formula
we have
$$ \langle c_1(J),[E_j]\rangle = E_j^2 - 2\genus(E_j) + 2 = E_j^2 + 2.$$

Now consider the set of tight contact structures $\xi_i$ on $Y$. Each $\xi_i$ is Stein fillable and
a Stein filling $(W^i,J^i)$ is given by taking a Legendrian surgery diagram, obtained from the plumbing tree, with the zigzags chosen in a certain way. Denote by $U^i_j$ the components of the corresponding Legendrian link and by $S^i_j$ the associated
surfaces in the Stein filling $(W^i,J^i)$. Notice that $W^i$ is diffeomorphic to $W$ by a diffeomorphism
which carries $S^i_j$ to $E_j$ for each $i$ and $j$. Also, since $E_j\cdot E_k$ is $0$ or $1$ if $j\neq k$, it follows that all components of each of the Legendrian links must have the same orientation.

Now, using the well-known identities
$$ (S^i_j)^2 = \tb(U^i_j)-1,\qquad \langle c_1(J^i), [S^i_j]\rangle = \rot(U^i_j), $$
observe that $\langle c_1(J^i),[S^i_j]\rangle = (S^i_j)^2 + 2$ precisely when $\rot(U^i_j) = \tb(U^i_j)+1$. Since the latter equality holds exactly when all the cusps of $U^i_j$ except one are up-cusps, it follows that
$\langle c_1(J),[E_j]\rangle = \langle c_1(J^i),[S^i_j]\rangle$
for each $j$ precisely when all the extra zigzags are chosen so that the additional cusps are all up-cusps, that is, when all the extra zigzags are
chosen on the same fixed side (which is determined by the orientation of the Legendrian unknots). The proof is now completed by appealing to the classification of contact structures on small Seifert-fibred $3$-manifolds
which are diffeomorphic to rational surface singularity links (cf. \cite{gh, wu}).
\end{proof}

\section{Final Remarks}
1) By Grauert's theorem, the small Seifert fibred $3$-manifold $Y = Y(e_0 ; r_1, r_2, r_3)$ is Milnor fillable if and only if the plumbing graph of $Y$ is negative definite. For example,
$Y(-2 ; \frac 12, \frac 23, \frac 56)$ is not Milnor fillable.

2) A necessary condition for the $3$-manifold $Y = Y(e_0 ; r_1, r_2, r_3)$ to be the link of a rational singularity is that $e_0\leq -2$; however, it is not sufficient. For example, $Y(-2 ; \frac 12, \frac 23, \frac {9}{11})$ is the link of an elliptic singularity.

\begin{figure}
\begin{center}
\includegraphics[scale=1.0]{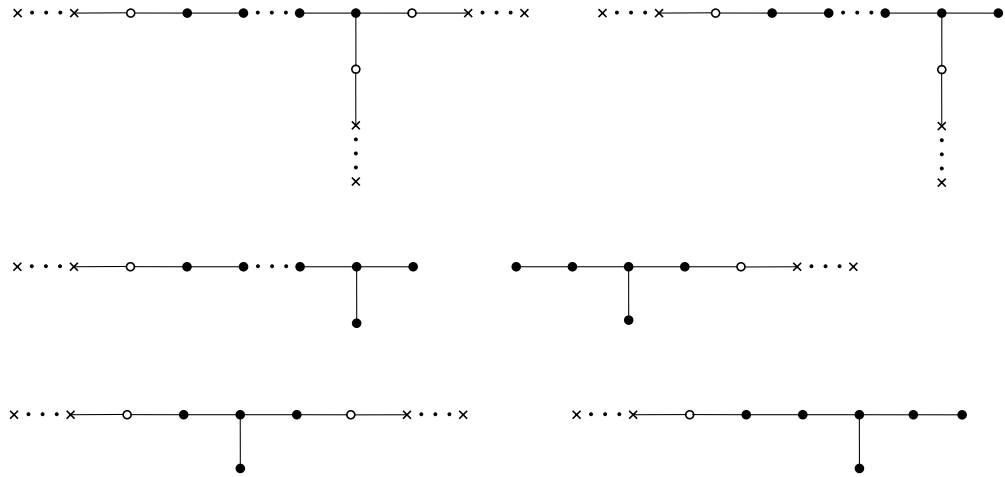}
\end{center}
\caption{Plumbing graphs of Seifert fibred $3$-manifolds whose canonical contact structures are supported by elliptic Milnor open books.}
\label{ellip}
\end{figure}

3) The canonical contact structure on $Y = Y(e_0 ; r_1, r_2, r_3)$ is supported by a planar Milnor open book if and only if $e_0\leq -3$. It is supported by an elliptic Milnor open book (and not a planar Milnor open book) if and only if the plumbing graph of $Y$ has one of the forms given in Figure \ref{ellip} or corresponds to one of the negative definite graphs $E_6, E_7$ or $E_8$.
Here the vertices are weighted as follows: vertices labelled $\bullet$ have weight $-2$, vertices labelled $\circ$ have weight less than $-2$ and vertices labelled $\times$ have weight less than $-1$. Also, the strings of vertices labelled $\times$ may be empty and further the leftmost $\circ$ may be absent in the first three families of graphs.

\begin{figure}
\begin{center}
\includegraphics[scale=1.0]{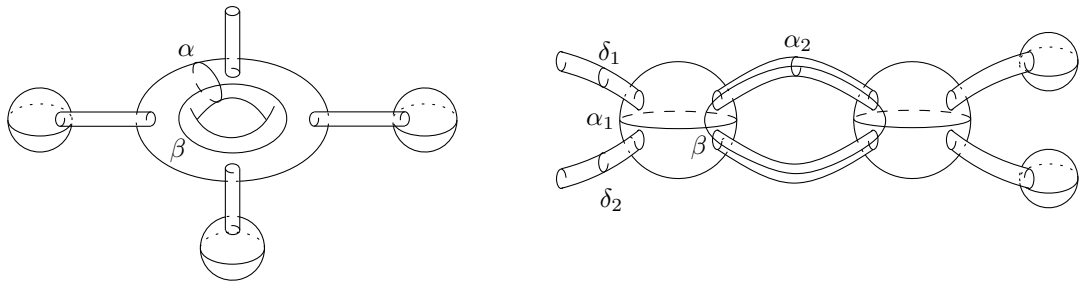}
\end{center}
\caption{Pages of two different Milnor open books on $Y(-2 ; \frac 12, \frac 12, \frac 12)$}
\label{D_4page}
\end{figure}

4) The minimal Milnor open book $\OB_{min}$ on $Y = Y(e_0 ; r_1, r_2, r_3)$ realizes $\Mg(\xi_{can})$, $\Mb(\xi_{can})$ and $\Mn(\xi_{can})$. In fact, it follows from the proof of Theorem \ref{Milmin} given in \cite{ab} that $\OB_{min}$ is the unique Milnor open book that realizes $\Mg(\xi_{can})$, $\Mb(\xi_{can})$ and $\Mn(\xi_{can})$.
Thus any other Milnor open book on $Y$ that realizes $\Mg(\xi_{can})$ cannot realize $\Mb(\xi_{can})$ and $\Mn(\xi_{can})$. For example, consider the $3$-manifold $Y=Y(-2 ; \frac 12, \frac 12, \frac 12)$, which is the link of the singularity $D_4$. The pages of two Milnor open books on $Y$ are given in Figure \ref{D_4page}. The first one is the minimal Milnor open book $\OB_{min}=\OB((1,2,1,1))$ with page a once-punctured torus and  monodromy $\phi=(\a\b)^3$; the second one is the Milnor open book $\OB=\OB((2,2,1,1))$ with page a twice-punctured torus and monodromy $\psi$ satisfying $\psi^2=\d_1\d_2\a_2^2$. Using the uniqueness result from \cite{bp} and the two-holed torus relation one can check that $\psi=\a_1\a_2\b\a_2^2\b\a_2$. It is easy to see that $\OB$ is related to $\OB_{min}$ by a single positive stabilization.

5) There are Milnor fillable contact $3$-manifolds such that  $$\sg(\xi_{can}) =\Mg(\xi_{can}), \sn(\xi_{can}) =
\Mn(\xi_{can}), \mbox{and} \bn(\xi_{can}) =\Mb(\xi_{can}).$$ For instance, the unique tight contact structure on the link of the singularity $E_8$ (with negative definite intersection form)  satisfies all the equalities above.

6) The invariants $\sg(\xi)$, $\sn(\xi)$ and $\bn(\xi)$ are independent in general for a contact structure $\xi$ as illustrated in \cite{be} (and \cite{el})
although
$$\Mn(\xi_{can})= 2 \Mg(\xi_{can}) -2 + \Mb(\xi_{can}).$$

7) There are examples of canonical contact structures such that Milnor genus, Milnor norm and Milnor binding number are arbitrarily large as we showed in Proposition~\ref{nobad} and  Theorem~\ref{arbit}.

8) There is a \emph{topological} lower bound for the support norm of an arbitrary contact structure on  a closed $3$-manifold $Y$.  Namely, the inequality  $\Hg(Y)-1  \leq \sn(\xi) $ holds for any contact structure $\xi$ on $Y$, where  $\Hg(Y)$ denotes the Heegaard genus of $Y$.   To see this,  we simply observe that every open book decomposition of a $3$-manifold  induces a Heegaard splitting. Moreover if a page of this open book  is a genus $g$ surface with $r$ boundary components then the Heegaard surface obtained by gluing two pages along their common boundary is a closed Riemann surface of genus $2g+r-1$.  It immediately follows that  $\Hg(Y) \leq 2g+r-1$ proving our claim  since  the negative of the Euler characteristic of a page of an open book  as above is equal to $2g+r-2$.

Moreover, it is well-known that the rank of $\pi_1(Y)$ is less than or equal to $\Hg(Y)$.  Hence we have  $$\rk (\pi_1(Y))  \leq \Hg(Y) \leq 1+ \sn(\xi)$$ for any contact structure $\xi$ on $Y$.

9) While there are no known example of a contact structure whose support genus is greater than one, there exist tight contact structures
with arbitrarily large support norm.  Let $Y_g$ denote the circle bundle over a closed genus $g$ surface whose Euler number is equal to $-1$ for $g \geq 0$.  There is a Stein fillable contact structure $\xi_g$ on $Y_g$ which  is supported by a horizontal open book with connected binding whose page has genus $g$ (see \cite{eo}).  Since $\Hg(Y_g)=2g$ (see \cite{boz}), we conclude that $\sn(\xi_g) \geq 2g-1$. The calculation of the
numbers $\sg(\xi_g)$ and $\bn(\xi_g)$ remains as a challenge.  Note that either $ \sg(\xi_g)$ or  $ \bn(\xi_g)$ is unbounded as $g \to \infty$, since $ 2g+1 \leq \sn(\xi_g) + 2  \leq 2 \sg(\xi_g) + \bn(\xi_g) $.

One can easily construct examples of  tight contact structures
with arbitrarily large binding numbers as we illustrate here.
Let $S_k = \#_k S^1 \times S^2$, for some positive integer $k$.
Consider the unique tight  contact structure  $\xi_k$ on $S_k$. It is known that
there is a planar open book compatible
with $\xi_k$ with $k+1$ binding components, where the monodromy is the identity map. We claim that this open book is minimal  in the sense that $\sg(\xi_k)=0$, $\bn(\xi_k) = k+1$, and $\sn(\xi_k) = k-1$.  Since Heegaard genus is additive under connected sum,   we have   $\Hg(S_k) = k$.  It follows that $\bn(\xi_k) - 1  \geq k$ and $\sn(\xi_k) \geq k-1$, since $\sg(\xi_k)=0$.  As  a consequence we conclude that $\sg(\xi_k)=0$, $\bn(\xi_k) = k+1$, and $\sn(\xi_k) = k-1$, since these numbers are realized by the aforementioned planar open book supporting $\xi_k$.

It would be interesting to find a sequence of tight contact structures on a fixed $3$-manifold with strictly increasing binding number and/or support norm.

\v \noindent {\bf {Acknowledgement}}: We would like to thank Andr\'{a}s  I. Stipsicz and John B. Etnyre
for suggesting some improvements after reading a preliminary version of this paper. The second author would also like to thank Vera V\'{e}rtesi and Ian Agol for helpful communications and the Mathematical Sciences Research Institute for its hospitality during the \textit{Symplectic and Contact Geometry and Topology} program 2009/10.
B.O. was partially supported by the TUBITAK-BIDEP-2219 research grant of the Scientific and Technological
Research Council of Turkey and the Marie Curie International Outgoing Fellowship 236639.

\end{document}